\documentclass[12pt,twoside]{article}
\usepackage{amsfonts}
\usepackage{ecit08}
\usepackage{amssymb}
\usepackage{latexsym}
\usepackage{graphicx}
\usepackage{multicol}
\usepackage{multirow}
\usepackage{color}
\usepackage{pstricks,pst-plot,pst-text,pst-tree,pst-eps,pst-fill,pst-node,pst-math}

\newcommand{\dbl}{[\![}
\newcommand{\dbr}{]\!]}
\newcounter{numerothm}
\newcounter{numerodef}
\newcounter{numerostat}
\newcounter{numeroprop}
\newcounter{numerocor}

\startpage{1}

\author[Ren\'e Lozi, Clarisse Fiol]{R. Lozi and C. Fiol} 

\address{\rm Ren\'e Lozi\\ \it
Laboratoire J.A. Dieudonn\'e - UMR du CNRS n$^{o}$6621\\
Universit\'e de Nice Sophia-Antipolis \\ Parc Valrose \\  06108 NICE CEDEX 2 \\FRANCE\\[.2in]
\rm Clarisse Fiol, Ren\'e Lozi \\ \it
IUFM C\'elestin Freinet - Universit\'e de Nice Sophia-Antipolis\\ 89, av. George V\\
06046 NICE CEDEX 1 \\FRANCE \\ \ \\ 
{\rm e-mail:} fiol@unice.fr \\ \hspace*{1.2cm} rlozi@unice.fr}

\title[Global Orbit Patterns for One Dimensional Dynamical Systems]{Global
Orbit Patterns for One Dimensional Dynamical Systems}

\date{}
\startpage{112}
\begin{document}
\maketitle \makeaddress \msc{37E15, 37M10, 70G60, 65P120}

\begin{abstract}
In this article, we study the behaviour of discrete one-dimensional
dynamical systems associated to functions on finite sets. We formalise the
global orbit pattern formed by all the periodic orbits (gop) as the ordered
set of periods when the initial value thumbs the finite set in the
increasing order. We are able to predict, using closed formulas, the number
of gop for the set $\mathcal{F}_N$ of all the functions on $X$. We also
explore by computer experiments special subsets of $\mathcal{F}_N$ in which
interesting patterns of gop are found.
\end{abstract}

\section{Introduction}\label{section1}

Dynamical systems $x_{n+1}=f_{\lambda}(x_n)$ associated to a function $%
f_{\lambda} : [a,b]\rightarrow [a,b]$ with $(a,b)\in \mathbb{R}^2$ on a real
interval are very well-known. Their periodic orbits, if any, $%
x_{n+1}^p=f_{\lambda}^p(x_n)$ are an important feature. There is a growing
interest in numerical analysis and industrial mathematics to study such
systems, or systems in higher dimensions \cite{Lozi1}. Very often, dynamical
systems in several dimensions are obtained coupling one-dimensional ones and
their properties are strongly linked.

Simple dynamical systems often involve periodic motion. Quasi-periodic or
chaotic motion is frequently present in more complicated dynamical systems.
The most famous theorem in this field of research is the Sharkovskii's
theorem, which addresses the existence of periodic orbits of continuous maps
of the real line into itself. This theorem was once proved toward the year
1962 and published only two years after, A.N. Sharkovskii \cite{Shark1}.

Mathematical results concerning periodic orbits are often obtained using
sets of real numbers. However, most of the time, as the complex behaviour of
chaotic dynamical systems is not explicitly tractable, mathematicians have
recourse to computer simulations. The main question which arises then is:
does these numerical computations are reliable?

As an example O.E. Lanford III \cite{Lanf} reports the results of some
computer experiments on the orbit structure of the discrete maps on a finite
set which arise when an expanding map of the circle is iterated "naively" on
the computer. \newline
Due to the discrete nature of floating points used by computers, there is a
huge gap between these results and the theoretical results obtained when
this map is considered on a real interval. (This gap can be narrowed in some
sense (i.e. avoiding the collapse of periodic orbits) in higher dimensions
when ultra weak coupling is used \cite{Lozi2, Lozi3}.)

Nowadays the claim is to understand precisely which periodic orbit can be
observed numerically in such systems. In a first attempt we study in this
paper the orbits generated by the iterations of a one-dimensional system on
a finite set $X_N$ with a cardinal $N$. The final goal of a good
understanding of the actual behaviour of dynamical systems acting on
floating numbers (i.e. the numbers used by computers) will be only reached
when this first step will be achieved.

On finite set, only periodic orbits can exist. For a given function we can
compute all the orbits, all together they form a global orbit pattern. We
formalise such a gop as the ordered set of periods when the initial value
thumbs the finite set in the increasing order. We are able to predict, using
closed formulas, the number of gop for the set $\mathcal{F}_N$ of all the
functions on $X$. We also explore by computer experiments special subsets of
$\mathcal{F}_N$ in which interesting patterns of gop are found.

This article is organized as follows: in section \ref{section2} we display
some examples already known of such computational discrepancies for the
logistic and circle maps. In section \ref{section3} we introduce a new
mathematical tool: the global orbit pattern, in order to describe more
precisely the behaviour of dynamical systems on finite sets. In section \ref%
{section4} we give some closed formulas giving the cardinal of classes of
gop of $\mathcal{F}_N$. In section \ref{section5} we study the case of sets
of functions with local properties versus their gop, in order to show the
significance of these new tools.

\section{Computational discrepancies}\label{section2}

\subsection{Approximated logistic map}

As an example of such collapsing effects, O.E. Lanford III, presents the
results of a sampling study in double precision of a discretization of the
logistic map
\begin{equation}  \label{eq0}
x_{n+1}=1-2x_n^2
\end{equation}

\noindent which has excellent ergodic properties. The precise discretization
studied is obtained by first exploiting evenness to fold the interval $[-1,0]
$ to $[0,1]$, i.e replacing (\ref{eq0}) by

\begin{equation}  \label{eq1}
x_{n+1}=|1-2x_n^2|
\end{equation}

On $[0, 1]$ it is not difficult to see that the folded map has the same
periods as the original one. The working interval is then shifted from $[0,
1]$ to $[1, 2]$ by translation, in order to avoid perturbation in numerical
experiment, by the special value 0. Then the translated folded map is
programmed in straightforward way. \newline
Out of 1,000 randomly chosen initial points,

\begin{itemize}
\item 890, i.e., the overwhelming majority, converged to the fixed point
corresponding to the fixed point-1 in the original representation (\ref{eq0}),

\item 108 converged to a cycle of period 3,490,273,

\item the remaining 2 converged to a cycle of period 1,107,319.
\end{itemize}

Thus, in this case at least, the very long-term behaviour of numerical
orbits is, for a substantial fraction of initial points, in flagrant
disagreement with the true behaviour of typical orbits of the original
smooth logistic map.

\subsection{Circle approximated maps}

In the same paper, O.E. Lanford III, studies very carefully the numerical
approximation of the map

\begin{equation}  \label{eq2}
x_{n+1}=2x_n+0.5x_n(1-x_n)~~~(\mathrm{mod}~1)~~0\leq x \leq 1
\end{equation}

\begin{figure}[tbp]
\begin{center}
\includegraphics*[width=6cm]{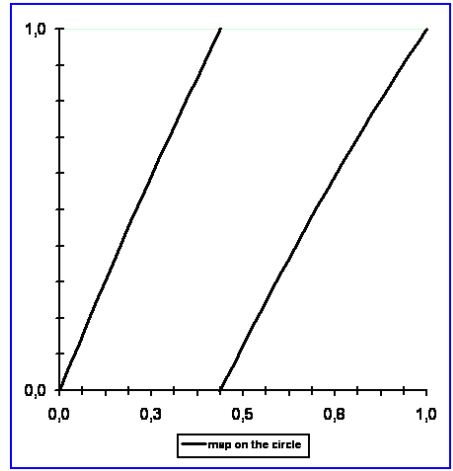}
\end{center}
\caption{Graph of the map $f(x)=2x+0.5x(1-x)$ (mod 1)} \label{fig1}
\end{figure}

It is perhaps better to think this map (see Figure \ref{fig1}) as acting on
the unit interval with endpoints identified, i.e., on the circle. Note that 
$f ^{\prime}(x)\geq 1.5$ everywhere, so $f$ is strictly expanding in a
particularly clean and simple sense. As a consequence of expansivity, this
mapping has about the best imaginable ergodic properties :

\begin{itemize}
\item it admits a unique invariant measure $\mu$ equivalent to Lebesgue
measure,

\item the abstract dynamical system $(f , \mu)$ is ergodic and in fact
isomorphic to a Bernoulli shift,

\item a central limit theorem holds.
\end{itemize}

One consequence of ergodicity of $f$ relative to $\mu$ is that, for Lebesgue
almost initial points $x$ in the unit interval, the corresponding orbit 
$f^{(n)}(x)$ is asymptotically distributed over the unit interval according
to $\mu$.

In numerical experiments performed by the author, the computer working with
fixed finite precision is able to represent finitely many points in the
interval in question. It is probably good, for purposes of orientation, to
think of the case where the representable points are uniformly spaced in the
interval. The true smooth map is then \emph{approximated} by a discretized
map, sending the finite set of representable points in the interval to
itself.

Describing the discretized mapping exactly is usually complicated, but it is
\emph{roughly} the mapping obtained by applying the exact smooth mapping to
each of the discrete representable points and "rounding" the result to the
nearest representable point.

O.E. Lanford III has done two kinds of experiments, in the first one
he uses double precision floating points, in the second uniformly
spaced points in the
interval with several order of discretization (ranging from $2^{22}$ to 
$2^{25}$) are involved. In each experiment the questions addressed are:

\begin{itemize}
\item how many periodic cycles are there and what are their periods ?

\item how large are their respective basins of attraction, i.e., for each
periodic cycle, how many initial points give orbits with eventually land on
the cycle in question?
\end{itemize}

\begin{table}[tbp]
\begin{center}
\begin{tabular}{rrr}
\hline
\multicolumn{3}{c}{Double precision (sampling)} \\
\multicolumn{3}{c}{$7$ cycles found} \\ \hline
period & basin size & relative size \\ \hline
27,627,856 & 517 & 51.7 \% \\
88,201,822 & 320 & 32.0 \% \\
4,206,988 & 147 & 14.7 \% \\
4,837,566 & 17 & 1.7 \% \\
802,279 & 8 & 0.8 \% \\
6,945,337 & 6 & 0.6 \% \\
2,808;977 & 1 & 0.1 \% \\ \hline
\end{tabular}%
\end{center}
\caption{Coexisting periodic orbits of mapping (\protect\ref{eq3}) for
double precision numbers.}
\label{tab:a14}
\end{table}

On one hand, for the experiments using ordinary (IEEE-754) double precision
- so that the working interval contains of the order of $10^{16}$
representable points - 1,000 initial points at random are used in order to
sample the orbit structure, determining the cycles to which they converge.

The computations were accomplished by shifting the working interval from the
initial interval [0, 1] to [1, 2] by translation, for the same previous
reason as for the logistic map, i.e., the map actually iterated was

\begin{equation}  \label{eq3}
x_{n+1}=2x_n+0.5x_n(x_n-1)(2-x_n)~~~(\mathrm{mod}~1)~~1\leq x \leq 2
\end{equation}

Some results are displayed in Table \ref{tab:a14}.\newline

On another hand, for relatively coarse discretizations the orbit
structure is determined completely, i.e., all the periodic cycles
and the exact sizes of their basins of attraction are found. Some
representative results are given in Tables \ref{tab:a0} to
\ref{tab:a00}. In theses tables, $N$ denotes the order of the
discretization, i.e., the representable points are the numbers,
$\frac{j}{N}$, with $0 \leq j < N$.

It seems that the computation of numerical approximations of the periodic
orbits leads to unpredictable results.

Many more examples could be given, but those given may serve to illustrate
the intriguing character of the results: the outcomes prove to be extremely
sensitive to the details of the experiment, but the results all have a
similar flavour: a relatively small number of cycles attracts near all
orbits, and the lengths of these significant cycles are much larger than one
but much smaller than the number of representable points. O. E. Lanford III
wrote that there are here regularities which ought to be understood.

\begin{table}[tbp]
\begin{center}
\begin{tabular}{rrr}
\hline
\multicolumn{3}{c}{$N=2^{23}=8,388,608$} \\
\multicolumn{3}{c}{$7$ cycles} \\ \hline
period & basin size & relative size \\ \hline
4,898 & 5,441,432 & 64.85 \% \\
1,746 & 2,946,734 & 35.13 \% \\
13 & 205 & $24.10^{-4} $ \% \\
6 & 132 & $16.10^{-4}$ \% \\
30 & 96 & $11.10^{-4}$ \% \\
4 & 8 & $< 1.10^{-6}$ \% \\
1 & 1 & $< 1.10^{-6}$ \% \\ \hline
\end{tabular}%
\end{center} \caption{Coexisting periodic orbits of mapping
(\protect\ref{eq3}) for the discretization $N=2^{23}$.}
\label{tab:a0}
\end{table}

\begin{table}[tbp]
\begin{center}
\begin{tabular}{rrr}
\hline
\multicolumn{3}{c}{$N=2^{24}=16,777,216$} \\
\multicolumn{3}{c}{$2$ cycles} \\ \hline
period & basin size & relative size \\ \hline
5,300 & 16,777,214 & 100 \% \\
1 & 2 & $< 1.10^{-6}$ \% \\ \hline
\end{tabular}%
\end{center}
\caption{Coexisting periodic orbits of mapping (\protect\ref{eq3})
for the discretization $N=2^{24}$.} \label{tab:a11}
\end{table}

\begin{table}[tbp]
\begin{center}
\begin{tabular}{rrr}
\hline
\multicolumn{3}{c}{$N=2^{24}-1=16,777,215$} \\
\multicolumn{3}{c}{$10$ cycles} \\ \hline
period & basin size & relative size \\ \hline
3,081 & 7,502,907 & 44.72 \% \\
699 & 3,047,369 & 18.16 \% \\
3,469 & 2,905,844 & 17.32 \% \\
1,012 & 2,774,926 & 16.54 \% \\
563 & 290,733 & 11.73 \% \\
2,159 & 221,294 & 1.32 \% \\
138 & 21,610 & 0.13 \% \\
421 & 12,477 & 0.07 \% \\
9 & 54 & $< 1.10^{-3}$ \% \\
1 & 1 & $< 1.10^{-7}$ \% \\ \hline
\end{tabular}
\end{center}
\caption{Coexisting periodic orbits of mapping (\protect\ref{eq3})
for the discretization $N=2^{24}-1$.} \label{tab:a12}
\end{table}

\begin{table}[tbp]
\begin{center}
\begin{tabular}{rrr}
\hline
\multicolumn{3}{c}{$N=2^{25}=33,554,432$} \\
\multicolumn{3}{c}{$8$ cycles} \\ \hline
period & basin size & relative size \\ \hline
4,094 & 32,114,650 & 95.71 \% \\
621 & 918,519 & 2.74 \% \\
283 & 516,985 & 1.54 \% \\
126 & 2,937 & $< 0.01$ \% \\
6 & 887 & $< 0.01$ \% \\
55 & 433 & $< 0.01$ \% \% \\
4 & 20 & $< 1.10^{-6}$ \% \\
1 & 1 & $< 1.10^{-6}$ \% \\ \hline
\end{tabular}%
\end{center}
\caption{Coexisting periodic orbits of mapping (\protect\ref{eq3})
for the discretization $N=2^{25}$.} \label{tab:a00}
\end{table}

In \cite{Dia}, P. Diamond and A. Pokrovskii, suggest that statistical
properties of the phenomenon of computational collapse of discretized
chaotic mapping can be modeled by random mappings with an absorbing centre.
The model gives results which are very much in line with computational
experiments and there appears to be a type of universality summarised by an
Arcsine law. The effects are discussed with special reference to the family
of mappings

\begin{equation}  \label{eq4}
x_{n+1}=1-|1-2x_n|^{\ell}~~~~0\leq x \leq 1 ~~~ 1 \leq \ell \leq 2
\end{equation}

Computer experiments show close agreement with prediction of the model.

However these results are of statistical nature, they do not give
accurate information on the exact nature of the orbits (e.g. length
of the shortest one, of the greater one, size of their basin of
attraction). It is why we consider the problem of computational
discrepancies in an original way in the next section.

\section{Global orbit pattern}

\label{section3}

We introduce in this section a new mathematical tool: the global orbit
pattern of a function, in order to describe more precisely which kind of
behaviour occurs in discretized dynamical systems on finite sets.

\subsection{Definitions}

Let $f$ be a map from $X$ onto $X$.\newline
We denote by $\forall x \in X$,~ $f^0(x)=x$ and $\forall p \geq 1,\forall x
\in X,~~f^p(x)=\underbrace{f \circ f \circ \ldots \circ f }_{\mathit{p}~~%
\textrm{\scriptsize{times}}}(x)$.\newline
Given $x_0 \in X$, we define the sequence $\lbrace x_i \rbrace$  of the
associated dynamical system by $x_{i+1}=f(x_i)$ for $i \geq 0$. Thus $%
x_i=f^i(x_0)$.\newline
The orbit of $x_0$ under $f$ is the set of points  $\mathcal{O}%
(x_0,f)=\lbrace f^i(x_0), i \geq 0 \rbrace$.\newline
The starting point $x_0$ for the orbit is called the initial value of the
orbit.\newline
A point $x$ is a fixed point of the map $f$ if $f(x)=x$.\newline
A point $x$ is a periodic point with period $p$ if $f^p(x)=x$ and $%
f^k(x)\neq x$ for all $k$ such that $0 \leq k <p$, $p$ is called the order
of $x$.\newline
If $x$ is periodic of order $p$, then the orbit of $x$ under $f$ is the
finite set \newline
$\lbrace x, f(x), f^2(x), \ldots, f^{p-1}(x) \rbrace$. We will call this set
the periodic orbit of order $p$ or a $p$-cycle.\newline
A fixed point is then a 1-cycle.\newline
The point $x$ is an eventually periodic point of $f$ with order $p$ if there
exists $K>0$ such that $\forall k\geq K$ $f^{k+p}(x)=f^k(x)$.\newline
$\forall x \in X$, we denote $\omega(x,f)$ the order of $x$ under $f$ or
simply $\omega(x)$ when the map $f$ involved is obvious.\newline
A subset $T$ of $X$ is invariant under $f$ if $f^{-1}(T)=T$. That is
equivalent to say that $T$ is invariant under $f$ if and only if $f(T)
\subset T$ and $f^{-1}(T)\subset T$.

\subsection{Map on finite set}

Along this paper, $N$ is a non-zero integer and $\sharp A$ stands for the
cardinal of any finite set $A$. In this article, we consider $X$ as an
ordered finite set with $N$ elements. We denote it $X_N$ and $X_N$ is
isomorphic to the interval in $\mathbb{N}$ $[\![ 0, N-1 ]\!] $. Then $\sharp
X_N=N$. Let $f$ be a map from $X_N$ into $X_N$. We denote by $\mathcal{F}_N$
the set of the maps from $X_N$ into $X_N$. Clearly, $\mathcal{F}_N$ is a
finite set and $\sharp \mathcal{F}_N= N^N$ elements. For all $x \in X_N$, $%
\mathcal{O}(x,f)$ is necessarily a finite set with at most $N$ elements.
Indeed, let us consider the sequence $\lbrace x, f(x), f^2(x), \ldots ,
f^{N-1}(x), f^N(x)\rbrace$ of the first $N+1$ iterated points. Thanks to the
Dirichlet's box principle, two elements are equals because $X_N$ has exactly
$N$ different values. Thus, every initial value of $X_N$ leads ultimately to
a repeating cycle. More precisely, if $x$ is a fixed point $\mathcal{O}(x,f)$
is the unique element $x$ and if $x$ is a periodic point with order $p$, $%
\mathcal{O}(x,f)$ has exactly $p$ elements. In this case, the orbit of $x$
under $f$ is the set $\mathcal{O}(x,f)= \lbrace x, f(x), f^2 (x), \ldots,
f^{p-1}(x) \rbrace$. If $x$ is an eventually periodic point with order $p$,
there exists $K>0$ such that $\forall k\geq K$ $f^{k+p}(x)=f^k(x)$. In this
case, the orbit of $x$ under $f$ is the set $\mathcal{O}(x,f)= \lbrace x,
f(x), f^2(x), \ldots, f^K(x), f^{K+1}(x), \ldots, f^{K+p-1}(x) \rbrace $.

\subsection{Equivalence classes}

\subsubsection{Components}

Let $f \in \mathcal{F}_N$. We consider on $X_N$ the relation $\sim$ defined
by : $\forall x, x^{\prime}\in X_N$,~ $x \sim x^{\prime}\Leftrightarrow
\exists k \in \mathbb{N}$ such that $f^k(x) \in \mathcal{O}(x^{\prime},f)$.
The relation $\sim$ is an equivalence relation on $X_N$. $\mathcal{S}_N/\sim$
is the collection of the equivalence classes that we will call components of
$X_N$ under $f$ which constitute a partition of $X_N$. The number of
components are given in \cite{Krus}. Asymptotic properties of the number of
cycles and components are studied in \cite{Mutaf}. For each component, we
take as class representative element the least element of the component. The
components will be written $T_N(x_0,f), \ldots, T_N(x_{p_{f,N}},f)$ where $%
x_i$ is the least element of $T_N(x_i,f)$. \newline

By analogy with real dynamical systems, we can define attractive and
repulsive components in discretized dynamical systems as follows.\newline

\stepcounter{numerodef} \textbf{{\emph{Definition}} \thenumerodef} ~~A
component is repulsive when it is a cycle. Otherwise, the component is
attractive.\newline

\textbf{\emph{Remark}} ~~In other words, a component is attractive when the
component contains at least an eventually periodic element. The
corresponding cycle is strictly contained in an attractive component.\newline

Examples are given in the tables \ref{tab:aa}, \ref{tab:bb} and \ref{tab:cc}.%
\newline

For instance, in the table \ref{tab:bb}, the function $f$ has
$\{8\}$ as orbit and $\{1,7,8\}$ as component which is attractive
because $\{1\}$ and $\{7\}$ are eventually periodic elements.
\newline

\begin{table}[!bth]
\begin{center}
\begin{tabular}{lc}
\hline
Function & orbit/component/nature \\ \hline
\begin{tabular}{ccc}
0 & $\rightarrow$ & 4 \\
1 & $\rightarrow$ & 1 \\
2 & $\rightarrow$ & 5 \\
3 & $\rightarrow$ & 4 \\
4 & $\rightarrow$ & 0 \\
5 & $\rightarrow$ & 9 \\
6 & $\rightarrow$ & 8 \\
7 & $\rightarrow$ & 5 \\
8 & $\rightarrow$ & 6 \\
9 & $\rightarrow$ & 7 \\
&  &
\end{tabular}
&
\begin{tabular}{ccc}
$\{0,4\}$ & $\{0,4,3\}$ & attractive \\
&  &  \\
$\{1\}$ & $\{1\}$ & repulsive \\
&  &  \\
$\{5,9,7\}$ & $\{2,5,9,7\}$ & attractive \\
&  &  \\
$\{6,8\}$ & $\{6,8\}$ & repulsive \\
&  &
\end{tabular}
\\ \hline
\end{tabular}%
\end{center}
\caption{Orbits and components of a function belonging to
$\mathcal{F}_{10}$ with gop $[2,1,3,2]_{10}$.} \label{tab:aa}
\end{table}

\begin{table}[!bth]
\begin{center}
\begin{tabular}{lc}
\hline
Function & orbit/component/nature \\ \hline
\begin{tabular}{ccc}
0 & $\rightarrow$ & 6 \\
1 & $\rightarrow$ & 8 \\
2 & $\rightarrow$ & 4 \\
3 & $\rightarrow$ & 9 \\
4 & $\rightarrow$ & 5 \\
5 & $\rightarrow$ & 2 \\
6 & $\rightarrow$ & 0 \\
7 & $\rightarrow$ & 1 \\
8 & $\rightarrow$ & 8 \\
9 & $\rightarrow$ & 3 \\
&  &
\end{tabular}
&
\begin{tabular}{ccc}
$\{0,6\}$ & $\{0,6\}$ & repulsive \\
&  &  \\
$\{8\}$ & $\{1,7,8\}$ & attractive \\
&  &  \\
$\{2,4,5\}$ & $\{2,4,5\}$ & repulsive \\
&  &  \\
$\{3,9\}$ & $\{3,9\}$ & repulsive \\
&  &
\end{tabular}
\\ \hline
\end{tabular}%
\caption{Orbits and components of another function belonging to $\mathcal{F}%
_{10}$ with gop $[2,1,3,2]_{10}$.}
\label{tab:bb}
\end{center}
\end{table}

\begin{table}[!bth]
\begin{center}
\begin{tabular}{lc}
\hline
Function & orbit/component/nature \\ \hline
\begin{tabular}{ccc}
0 & $\rightarrow$ & 9 \\
1 & $\rightarrow$ & 6 \\
2 & $\rightarrow$ & 9 \\
3 & $\rightarrow$ & 8 \\
4 & $\rightarrow$ & 3 \\
5 & $\rightarrow$ & 7 \\
6 & $\rightarrow$ & 6 \\
7 & $\rightarrow$ & 5 \\
8 & $\rightarrow$ & 4 \\
9 & $\rightarrow$ & 2 \\
&  &
\end{tabular}
&
\begin{tabular}{ccc}
$\{2,9\}$ & $\{0,2,9\}$ & attractive \\
&  &  \\
$\{6\}$ & $\{1,6\}$ & attractive \\
&  &  \\
$\{3,8,4\}$ & $\{3,8,4\}$ & repulsive \\
&  &  \\
$\{5,7\}$ & $\{5,7\}$ & repulsive \\
&  &
\end{tabular}
\\ \hline
\end{tabular}%
\end{center}
\caption{Orbits and components of a third function belonging to $\mathcal{F}%
_{10}$ with gop $[2,1,3,2]_{10}$.}
\label{tab:cc}
\end{table}

\subsubsection{Order of elements}

Here are some remarks on the order of elements of components.\newline

\textbf{\emph{Remark}} ~~The order of every element of a component is the
length of its inner cycle.\newline

\stepcounter{numerodef} \textbf{{\emph{Definition}} \thenumerodef} ~~For all
$x \in X_N$, there exists $i \in [\![ 0, p_{f,N} ]\!] $ such that $x$
belongs to the component $T_N(x_i,f)$. Then $\omega(x,f)$ is equal to the
order $\omega(x_i,f)$. \newline

\textbf{\emph{Remark}} ~~For all $i \in [\![ 0,p_{f,N} ]\!]$, $T_N(x_i,f)$
is an invariant subset of $X_N$ under $f$.\newline

In the example given in the table \ref{tab:aa}, the order of the
element $\{5\}$ is 3, the order of the element $\{1\}$ is $1$, the
order of the element $\{8\}$ is 2. The elements $\{2\}$ and $\{5\}$
have the same order.

\subsection{Definition of global orbit pattern}

For each $f \in \mathcal{F}_N$, we can determine the components of $X_N$
under $f$. For each component, we determine the order of any element. Thus,
for each $f \in \mathcal{F}_N$, we have a set of orders that we will denote $%
\Omega(f, N)$. Be given $f$, there exist $p_{f,N}$ components and $p_{f,N}$
representative elements such that $x_0 < x_1 < \ldots <x_{p_{f,N}}$.\newline

For each $f \in \mathcal{F}_N$, the sequence $[\omega(x_0), \omega(x_1),
\ldots, \omega(x_{p_{f,N}}); f]_{\mathcal{F}_N}$ with $x_0 < x_1 < \ldots
<x_{p_{f,N}}$ will design  the global orbit pattern of $f \in \mathcal{F}_N$%
. \newline

We will write $gop(f)=[\omega(x_0), \omega(x_1), \ldots,
\omega(x_{p_{f,N}}); f]_{\mathcal{F}_N}$.\newline

When the reference to $f \in \mathcal{F}_N$ is obvious, we will write
shortly \newline
$gop(f)=[\omega(x_0), \omega(x_1), \ldots,  \omega(x_{p_{f,N}})]_N$ or $%
gop(f)=[\omega(x_0), \omega(x_1), \ldots, \omega(x_{p_{f,N}})]$ .\newline

For example, the same gop associated to the functions given in the tables %
\ref{tab:aa}, \ref{tab:bb} and \ref{tab:cc} is $[2,1,3,2]_{10}$.\newline

Another example is given in the table \ref{tab:dd}. In that example, we have
$\omega(0)=2$, $\omega(3)=1$, $\omega(4)=4$. \newline
\newline

\begin{table}[!bth]
\begin{center}
\begin{tabular}{lc}
\hline
Function & orbit/component/nature \\ \hline
\begin{tabular}{ccc}
0 & $\rightarrow$ & 1 \\
1 & $\rightarrow$ & 0 \\
2 & $\rightarrow$ & 0 \\
3 & $\rightarrow$ & 3 \\
4 & $\rightarrow$ & 5 \\
5 & $\rightarrow$ & 6 \\
6 & $\rightarrow$ & 7 \\
7 & $\rightarrow$ & 4 \\
&  &
\end{tabular}
&
\begin{tabular}{ccc}
$\{0,1\}$ & $\{0,1,2\}$ & attractive \\
&  &  \\
$\{3\}$ & $\{3\}$ & repulsive \\
&  &  \\
$\{4,5,6,7\}$ & $\{4,5,6,7\}$ & repulsive \\
&  &  \\
&  &
\end{tabular}
\\ \hline
\end{tabular}%
\end{center}
\caption{Orbits and components of a function belonging to $\mathcal{F}_{8}$
with gop $[2,1,4]_{8}$.}
\label{tab:dd}
\end{table}

\stepcounter{numerodef} \textbf{{\emph{Definition}} \thenumerodef} ~~The set
of all the global orbit patterns of $\mathcal{F}_N$ is called $\mathcal{G}(%
\mathcal{F}_N)$.\newline

For example, for $N=4$, the set $\mathcal{G}(\mathcal{F}_4)$ is $\{ [1];
[1,1]; [1,1,1]; [1,2]; [1,1,1,1];$ $[1,1,2]; [1,2,1]; [1,3]; [2]; [2,1];
[2,1,1]; [2,2]; [3]; [3,1]; [4]\}$.

\subsection{Class of gop}

We give the following definitions : \newline

\stepcounter{numerodef} \textbf{{\emph{Definition}} \thenumerodef} ~~Let be $%
A=[\omega_1, \ldots ,\omega_p]_N$ a gop. Then the class of $A$, written $%
\overline{A}$, is the set of all the functions $f \in \mathcal{F}_N$ such
that the global orbit pattern associated to $f$ is $A$.\newline

For example, for $N=10$, the class of the gop $\overline{[2,1,3,2]}_{10}$
contains the following few of many functions defined in the tables \ref%
{tab:aa}, \ref{tab:bb} and \ref{tab:cc}. The periodic orbit which are
encountered have the same length nevertheless there are different.\newline

\stepcounter{numerodef} \textbf{{\emph{Definition}} \thenumerodef} ~~Let be $%
A=[\omega_1, \ldots ,\omega_p]_N$ a gop. \newline
Then the modulus of $A$ is $|A| = \sum\limits_{k=1}^p \omega_k$.\newline

\textbf{\emph{Remark}} ~~$\left|[\omega_1, \ldots ,\omega_p]_N\right| \leq N$%
.\newline

\subsection{Threshold functions}

\subsubsection{Ordering the discrete maps}

\stepcounter{numerothm} \textbf{{\emph{Theorem}} \thenumerothm} ~~The sets $%
\mathcal{F}_N$ and $[\![ 1,N^N ]\!]$ are isomorphic.\newline

\proof{ We define the function $\phi$ from $\mathcal{F}_N$ to $\dbl
1,N^N \dbr$ by : for each  $f \in \mathcal{F}_N$, $\phi(f)$ is the
integer $n$ such that $n = \sum\limits_{k=0}^{N-1} f(k)N^{N-1-k} + 1$.\\
Then
$\phi$ is well defined because $n \in \dbl 1,N^N \dbr$. \\
Let $n$ be a given integer between 1 and $N^N$. We convert $n-1$ in
base $N$ : there exists a unique $N$-tuple $(a_{n-1,0}; a_{n-1,1};
\ldots ; a_{n-1,N-1}) \in \dbl 0,N-1\dbr^N$ such that
$\overline{n-1}^N = \sum\limits_{i=0}^{N-1} a_{n-1,N-1-i}
N^{N-i-1}$. We can thus define the map $f_n$ with : $\forall i \in
X_N$, $f_n(i)=a_{n-1,N-i-1}$. Then $\phi$ is one to one. }\newline

\textbf{\emph{Remark}} ~~This implies $\mathcal{F}_N$ is totally ordered.%
\newline

\stepcounter{numerodef} \textbf{{\emph{Definition}} \thenumerodef} ~~Let $f
\in \mathcal{F}_N$. Then $n=\sum\limits_{k=0}^{N-1} f(k)N^{N-1-k}+ 1$ is
called the rank of $f$. \newline

\subsubsection{Threshold functions}

Be given a global orbit pattern $A$, we are exploring the class $\overline{A}
$.\newline

\stepcounter{numerothm} \textbf{{\emph{Theorem}} \thenumerothm} ~~For every $%
A \in \mathcal{G}(\mathcal{F}_N)$, the class $\overline{A}$ has a unique
function with minimal rank.\newline

\stepcounter{numerodef} \textbf{{\emph{Definition}} \thenumerodef} ~~For
every class $\overline{A} \in \mathcal{G}(\mathcal{F}_N)$, the function
defined by the previous theorem will be called the threshold function for
the class $\overline{A}$ and will be denoted by $Tr(A)$ or $Tr(\overline{A})$%
. \newline

To prove the theorem, we need the following definition : \newline

\stepcounter{numerodef} \textbf{{\emph{Definition}} \thenumerodef} ~~Let $f
\in \mathcal{F}_N$ be a function. Let $p$ be a non zero integer smaller than
$N$. Let be $x_1, \ldots, x_{p}$ $p$ consecutive elements of $X_N$. Then $%
x_1, \ldots, x_{p}$ is a canonical $p$-cycle in relation to $f$ if $\forall
j \in [\![ 1, p-1 ]\!]$ , $f(x_j)=x_{j+1}$ and $f(x_p)=x_1$.\newline

\proof{Let $[\omega_1, \ldots ,\omega_p]$ be a global orbit pattern
of $\mathcal{G}(\mathcal{F}_N)$. We construct a specific function
$f$ belonging to the class $\overline{[\omega_1, \ldots ,\omega_p]}$
and we prove that the function so obtained is the smallest with
respect to the order on $\mathcal{F}_N$. With the first $\omega_1$
elements of $\dbl 0,N-1 \dbr$, that is the set of integers $\dbl
0,\omega_1-1 \dbr$, we construct the canonical $\omega_1$-cycle: if
$\omega_1=1$, we define $f(0)=0$, else $f(0)=1$, $f(1)=2$, $\ldots$,
$f(\omega_1-2)=\omega_1-1$, $f(\omega_1-1)=0$.
\\Then $\forall j \in \dbl \omega_1-1, \omega_1+N-s-1 \dbr$, we define
$f(j)=0$. \\Then with the next $\omega_2$ integers $\dbl
\omega_1+N-s,\omega_1+N-s+\omega_2-1 \dbr$ we construct the
canonical $\omega_2$-cycle. We keep going on constructing for all
 $ j \in  \dbl 3,p \dbr$ the canonical $\omega_j$-cycle.\\
 In consequence, we have found a function $f$ belonging to the class 
 $\overline{[\omega_1, \ldots ,\omega_p]}$. \\
Assume there exists a function $g \in \mathcal{F}_N$ belonging to
the class of $f$ such that $g < f$. Let $I= \{i \in \dbl 0,N-1 \dbr
~\textrm{such that}~ f(i) \neq 0 \} $. As $g <f$, there exists $i_0
\in I$ such that $g(i_0) <f(i_0)$. There exists also $j_0$ such that
$i_0 \in \omega_{j_0}$. If $f(i_0)=i_0$, then $\omega_{j_0}=1$,
$g(i_0)<i_0$ and then $g(i_0) \notin \omega_{j_0}$. Then the global
orbit pattern of $g$ doesn't contain anymore 1 as cycle. The global
orbit pattern of $g$ is different from the global orbit pattern of
$f$. If $f(i_0)=i_0+1$, then $g(i_0)\leq i_0$. Either $g(i_0)=i_0$
and then the global orbit pattern of $g$ is changed, or $g(i_0)<i_0$
and we are in the same situation as previously. Thus, in any case,
the smallest function belonging to the class $[\overline{\omega_1,
\ldots ,\omega_p]}$ is the
one constructed in the first part of the proof.}\newline

The proof of the theorem gives an algorithm of construction of the threshold
function associated to a given gop.\newline
The threshold function associated to the gop $[2,1,3,2]_{10}$ is explained
in the table \ref{tab:ee}. Its rank is $n=1,000,467,598$.\newline

\begin{table}[tbp]
\begin{center}
\begin{tabular}{lllll}
\hline
First step & Second step & Third step & Fourth step & Fifth step \\ \hline
\begin{minipage}{2,1cm} Construction of the first canonical 2-cycle\\
\begin{tabular}{ccc} 0 & $\rightarrow$ & 1\\ 1 & $\rightarrow$ & 0\\ 2 &
$\rightarrow$ & \\ 3 & $\rightarrow$ & \\ 4 & $\rightarrow$ & \\ 5 &
$\rightarrow$ & \\ 6 & $\rightarrow$ & \\ 7 & $\rightarrow$ & \\ 8 &
$\rightarrow$ & \\ 9 & $\rightarrow$ & \\ \end{tabular} \end{minipage} & %
\begin{minipage}{2,1cm} Construction of the last canonical 2-cycle\\
\begin{tabular}{ccc} 0 & $\rightarrow$ & 1\\ 1 & $\rightarrow$ & 0\\ 2 &
$\rightarrow$ & \\ 3 & $\rightarrow$ & \\ 4 & $\rightarrow$ & \\ 5 &
$\rightarrow$ & \\ 6 & $\rightarrow$ & \\ 7 & $\rightarrow$ & \\ 8 &
$\rightarrow$ & 9\\ 9 & $\rightarrow$ & 8\\ \end{tabular} \end{minipage} & %
\begin{minipage}{2,1cm} Construction of the canonical 3-cycle\\
\begin{tabular}{ccc} 0 & $\rightarrow$ & 1\\ 1 & $\rightarrow$ & 0\\ 2 &
$\rightarrow$ & \\ 3 & $\rightarrow$ & \\ 4 & $\rightarrow$ & \\ 5 &
$\rightarrow$ & 6\\ 6 & $\rightarrow$ & 7\\ 7 & $\rightarrow$ & 5\\ 8 &
$\rightarrow$ & 9\\ 9 & $\rightarrow$ & 8\\ \end{tabular} \end{minipage} & %
\begin{minipage}{2,1cm} Construction of the canonical 1-cycle\\
\begin{tabular}{ccc} 0 & $\rightarrow$ & 1\\ 1 & $\rightarrow$ & 0\\ 2 &
$\rightarrow$ & \\ 3 & $\rightarrow$ & \\ 4 & $\rightarrow$ & 4\\ 5 &
$\rightarrow$ & 6\\ 6 & $\rightarrow$ & 7\\ 7 & $\rightarrow$ & 5\\ 8 &
$\rightarrow$ & 9\\ 9 & $\rightarrow$ & 8\\ \end{tabular} \end{minipage} & %
\begin{minipage}{2,1cm} Filling the remaining images with 0\\
\begin{tabular}{ccc} 0 & $\rightarrow$ & 1\\ 1 & $\rightarrow$ & 0\\ 2 &
$\rightarrow$ & 0\\ 3 & $\rightarrow$ & 0\\ 4 & $\rightarrow$ & 4\\ 5 &
$\rightarrow$ & 6\\ 6 & $\rightarrow$ & 7\\ 7 & $\rightarrow$ & 5\\ 8 &
$\rightarrow$ & 9\\ 9 & $\rightarrow$ & 8\\ \end{tabular} \end{minipage} \\
\hline
\end{tabular}%
\end{center}
\caption{Algorithm for the threshold function construction for the gop $%
[2,1,3,2]_{10}$.}
\label{tab:ee}
\end{table}

\stepcounter{numerothm} \textbf{{\emph{Theorem}} \thenumerothm} ~~There are
exactly $2^N-1$ different global orbit patterns in $\mathcal{F}_N$.\newline

That is $\sharp \mathcal{G}(\mathcal{F}_N)=2^N-1$.\newline

For example, for $N=4$, $\sharp \mathcal{G}(\mathcal{F}_4)=2^4-1=15$.\newline

\proof{ Let $p$ be an integer between 1 and $N$. Consider the set
$L(p,N)$ of $p$-tuples $(\alpha_1, \ldots, \alpha_p) \in
(\mathbb{N}^*)^p $ such that $\alpha_1+ \ldots+ \alpha_p \leq N$.
\\We write $L(N) = \lbrace L(p,N), p=1 \ldots N \rbrace$. $L(N)$ and
$\mathcal{G}(\mathcal{F}_N)$ have the same elements. Then \\
\begin{center}
$ \displaystyle \sharp \mathcal{G}(\mathcal{F}_N)= \sum
\limits_{p=1}^{p=N} \sharp L(p,N)= \sum \limits_{p=1}^{p=N}
\left(\begin{array}{c}
  N \\
  p \\
\end{array}\right)=2^N-1$.
\end{center}}

\subsubsection{Ordering the global orbit patterns}

We define an order relation on $\mathcal{G}(\mathcal{F}_N)$.\newline

\stepcounter{numeroprop} \textbf{{\emph{Proposition}} \thenumeroprop} ~~Let $%
A$ and $B$ be two global orbit patterns of $\mathcal{G}(\mathcal{F}_N)$.%
\newline
We define the relation $\prec$ on the set $\mathcal{G}(\mathcal{F}_N)$ by
\[
A \prec B~~ \mathit{iff} ~~Tr(A)<Tr(B)
\]
Then the set $(\mathcal{G}(\mathcal{F}_N),\prec)$ is totally ordered.\newline

\proof{As the order  $\prec$ refers to the natural order of
$\mathbb{N}$, the proof is obvious.}\newline

Let $r \geq 1$, $p\geq 1$ be two integers. Let $[\omega_1, \ldots ,\omega_p]$
and $[\omega^{\prime}_1, \ldots, \omega^{\prime}_r]$ be two global orbit
patterns of $\mathcal{G}(\mathcal{F}_N)$. For example, if $p<r$, in order to
compare them, we admit that we can fill $[\omega_1, \ldots ,\omega_p]$ with $%
r-p$ zeros and write $[\omega_1, \ldots ,\omega_p]=[\omega_1, \ldots
,\omega_p, 0, \ldots, 0]$.\newline

\stepcounter{numeroprop} \textbf{{\emph{Proposition}} \thenumeroprop} ~~Let $%
r \geq 1$, $p\geq 1$ be two integers such that $p \leq r$. Let $A
=[\omega_1, \ldots ,\omega_p]$ and $B =[\omega^{\prime}_1, \ldots,
\omega^{\prime}_r]$ be two global orbit patterns. \newline
\begin{minipage}{14cm}
\begin{itemize}
\item If $r=p=1$ and $\omega_1 < \omega'_1$ then $A \prec B$.
\item If $r \geq 2$ then
\begin{itemize}
  \item[$\ast$]  If $\omega_1 < \omega'_1$ then $A \prec B$.
    \item[$\ast$]  If $\omega_1 = \omega'_1$ then there exists $K \in \dbl 2;r \dbr$
such that $\omega_K \neq \omega'_K$ and $\forall i < K$, $\omega_i =
\omega'_i$.
\begin{itemize}
  \item[$\bullet$] If $|A| < |B|$, then
$A \prec B$.
   \item[$\bullet$] If $|A| = |B|$, then if
$\omega_K < \omega'_K$ then $A \prec B$.
\end{itemize}
\end{itemize}
\end{itemize}
\end{minipage} \newline
\newline

For example, for $N=4$, the global orbit patterns are in increasing order : $%
[1] \prec [1,1] \prec [1,1,1] \prec [1,2] \prec [1,1,1,1] \prec [1,1,2]
\prec [1,2,1] \prec [1,3] \prec [2] \prec [2,1] \prec [2,1,1] \prec [2,2]
\prec [3] \prec [3,1] \prec [4] $.

\subsubsection{Algorithm for ordering the global orbit patterns : a
pseudo-decimal order}

The table \ref{tab:a} gives a method for ordering the gop : indeed, we
consider each gop as if each one represents a decimal number: we begin to
order them in considering the first order $\omega_1$. Considering two gops $%
A= [\omega_1, ..., \omega_p]$ and $A^{\prime}=[\omega_1^{\prime}, ...,
\omega_r^{\prime}]$, if $\omega_1<\omega_1^{\prime}$, then $A \prec
A^{\prime}$. For example, $[1,3] \prec [2,1]$. If $\omega_1=\omega_1^{\prime}
$ and $|A|-\omega_1<|A^{\prime}|-\omega_1^{\prime}$, then $A \prec A^{\prime}
$. For example to compare the gop $[1,2]$ and the gop $[1,1,1,1]$, we say
that the first order $\omega_1$ stands for the unit digit - which is $%
\omega_1=1$ here, then the decimal digits are respectively $0.2$ and $0.111$%
. We calculate for each of them the modulus-$\omega_1$: we find $|[1,2]| -1
= 2$ and $|[1,1,1,1]| -1 = 3$, thus $[1,2] \prec [1,1,1,1]$. Finally, if $%
\omega_1=\omega_1^{\prime}$ and $|A|-\omega_1=|A^{\prime}|-\omega_1^{\prime}$%
, then also we use the order of the decimal part. For example, $[1,1,1]
\prec [1,2]$ because $1.11 < 1.2$. Applying this process, we have the
sequence of the ordered gop for $N=4$ given in the previous paragraph.%
\newline

\begin{table}[!bth]
\begin{center}
\begin{tabular}{lcc||lcc}
\hline
Gop & Modulus & Modulus-$\omega_1$ & Gop & Modulus & Modulus-$\omega_1$ \\
\hline
$[1] $ & 1 & 0 & $[2]$ & 2 & 0 \\
$[1,1]$ & 2 & 1 & $[2,1]$ & 3 & 1 \\
$[1,1,1]$ & 3 & 2 & $[2,1,1]$ & 4 & 2 \\
$[1,2]$ & 3 & 2 & $[2,2]$ & 4 & 2 \\
$[1,1,1,1]$ & 4 & 3 &  &  &  \\
$[1,1,2]$ & 4 & 3 & $[3]$ & 3 & 0 \\
$[1,2,1]$ & 4 & 3 & $[3,1]$ & 4 & 1 \\
$[1,3]$ & 4 & 3 & $[4]$ & 4 & 0 \\ \hline
\end{tabular}%
\end{center}
\caption{Ordered gop for $N=4$ with modulus and modulus-$\protect\omega_1$}
\label{tab:a}
\end{table}

\begin{figure}[tbp]
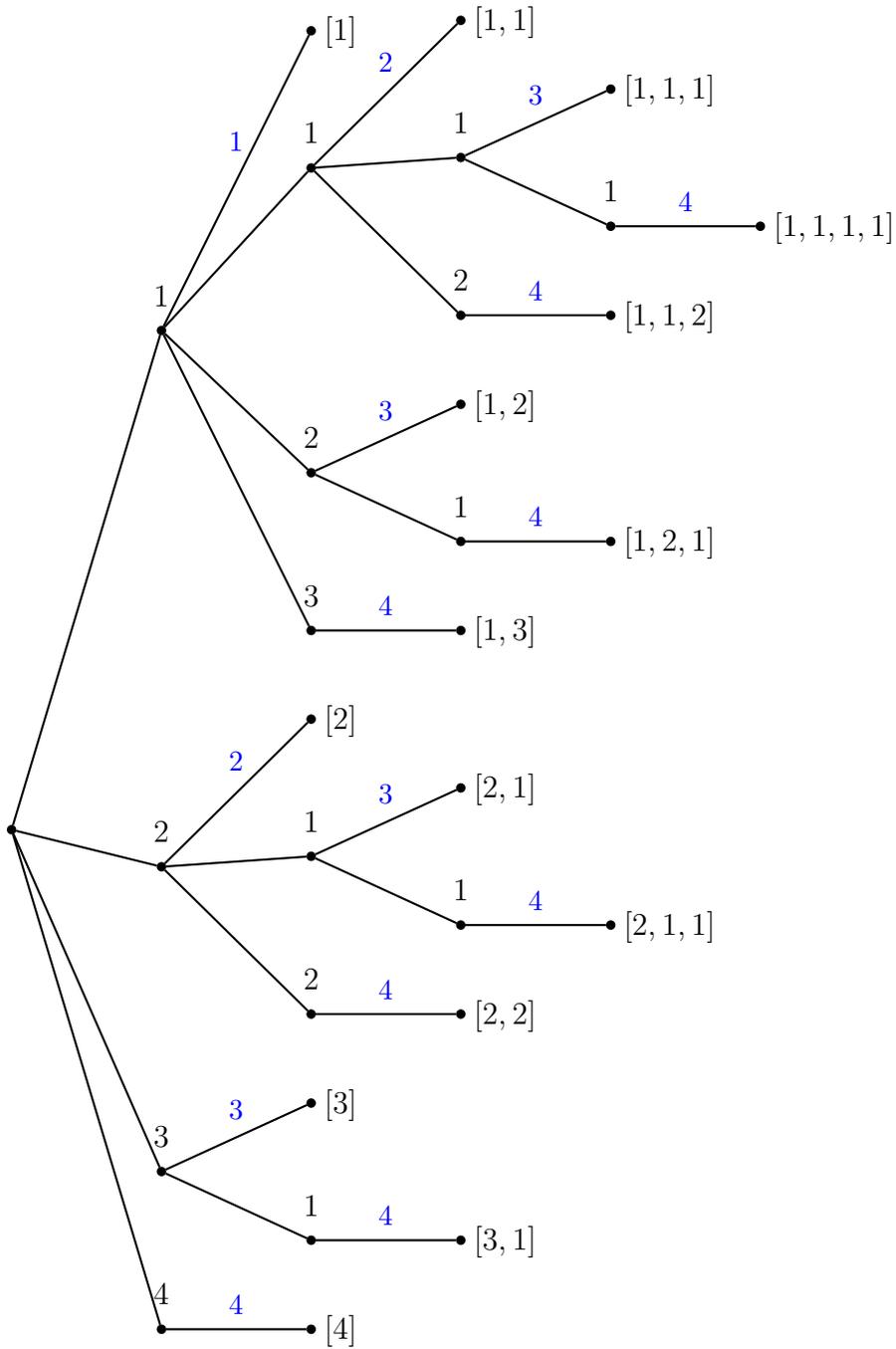

\psset{nodesep=0mm,levelsep=20mm,treesep=10mm}
\pstree[treemode=R]{\Tdot} { \pstree
{\Tdot~[tnpos=a]{$1$}\taput{\small $$}} {
\Tdot~[tnpos=r]{$[1]$}\taput{\small $\color{blue}{1}$} \pstree
{\Tdot~[tnpos=a]{$1$}\taput{\small $$}} {
\Tdot~[tnpos=r]{$[1,1]$}\taput{\small $\color{blue}{2}$} \pstree
{\Tdot~[tnpos=a]{$1$}\taput{\small $$}} {
\Tdot~[tnpos=r]{$[1,1,1]$}\taput{\small $\color{blue}{3}$} \pstree
{\Tdot~[tnpos=a]{$1$}\tbput{\small $$}} {
\Tdot~[tnpos=r]{$[1,1,1,1]$}\taput{\small $\color{blue}{4}$} } }
\pstree {\Tdot~[tnpos=a]{$2$}\tbput{\small $$}} {
\Tdot~[tnpos=r]{$[1,1,2]$}\taput{\small $\color{blue}{4}$} } }
\pstree {\Tdot~[tnpos=a]{$2$}\taput{\small $$}} {
\Tdot~[tnpos=r]{$[1,2]$}\taput{\small $\color{blue}{3}$} \pstree
{\Tdot~[tnpos=a]{$1$}\tbput{\small $$}} {
\Tdot~[tnpos=r]{$[1,2,1]$}\taput{\small $\color{blue}{4}$} } }
\pstree {\Tdot~[tnpos=a]{$3$}\taput{\small $$}} {
\Tdot~[tnpos=r]{$[1,3]$}\taput{\small $\color{blue}{4}$} } } \pstree
{\Tdot~[tnpos=a]{$2$}\taput{\small $$}} {
\Tdot~[tnpos=r]{$[2]$}\taput{\small $\color{blue}{2}$} \pstree
{\Tdot~[tnpos=a]{$1$}\taput{\small $$}} {
\Tdot~[tnpos=r]{$[2,1]$}\taput{\small $\color{blue}{3}$} \pstree
{\Tdot~[tnpos=a]{$1$}\tbput{\small $$}} {
\Tdot~[tnpos=r]{$[2,1,1]$}\taput{\small $\color{blue}{4}$} } }
\pstree {\Tdot~[tnpos=a]{$2$}\tbput{\small $$}} { \Tdot~[tnpos=r]{$
[2,2]$}\taput{\small $\color{blue}{4}$} } } \pstree
{\Tdot~[tnpos=a]{$3$}\taput{\small $$}} {
\Tdot~[tnpos=r]{$[3]$}\taput{\small $\color{blue}{3}$} \pstree
{\Tdot~[tnpos=a]{$1$}\tbput{\small $$}} {
\Tdot~[tnpos=r]{$[3,1]$}\taput{\small $\color{blue}{4}$} } } \pstree
{\Tdot~[tnpos=a]{$4$}\taput{\small $$}} {
\Tdot~[tnpos=r]{$[4]$}\taput{\small $\color{blue}{4}$} } }\newline
\caption{Tree for the construction of the order on $\mathcal{G}(\mathcal{F}_4)$}
\label{fig2}
\end{figure}

For example, for $N=4$, we construct a tree (Figure \ref{fig2}) : each
vertex is an ordered orbit, the modulus of the gop is written on the last
edge.\newline
However, the sequence of ordered gop differs from the natural downward
lecture of the tree and has to be done following the algorithm.

\section{Numbering the discrete maps}

\label{section4}

In this section we emphasize some closed formulas giving the cardinal of
classes of gop. Recalling first the already known formula for the class $[%
\overline{\underbrace{1, \ldots, 1}_{k ~\textrm{\scriptsize{times}}}}]_N$ for
which we give a detailed proof, we consider the case were the class
possesses exactly one $k$-cycle, the case with only two cycles belonging to
the class and finally the main general formula of any cycles with any
length. We give rigorous proof of all. The general formula is very
interesting in the sense that even using computer network it is impossible
to check every function of $\mathcal{F}_N$ when $N$ is larger than 100 or
may be 50.

\subsection{Discrete maps with 1-cycle only}

The theorem 4 gives the number of discrete maps of $\mathcal{F}_N$ which
have only fixed points and no cycles of length greater than one. This
formula is explicit in \cite{Knuth} and \cite{Purd}. A complete proof is
given here in detail.\newline

\label{th4}\stepcounter{numerothm} \textbf{{\emph{Theorem}} \thenumerothm}
~~ Let $k$ be an integer between 1 and $N$. The number of functions whose
global orbit pattern is $[\underbrace{1, \ldots, 1}_{k ~\textrm{\scriptsize{times}}%
}]_N$ (i.e. belonging to the class $[\overline{\underbrace{1, \ldots, 1}_
{k ~\textrm{\scriptsize{times}}}}]_N$) is $\left(%
\begin{array}{c}
N-1 \\
N-k \\
\end{array}
\right) N^{N-k}$.\newline
\newline

That is $\sharp [\overline{\underbrace{1, \ldots, 1}_{\mathit{k} ~\textrm{\scriptsize{times}}
}}]_N = \left(
\begin{array}{c}
N-1 \\
N-k \\
\end{array}
\right) N^{N-k}$.

\proof{ Let $k$ be a non-zero integer. Let $f$ be a function of
$\mathcal{F}_N$. There are $\left(\begin{array}{c}
  N \\
  k \\
\end{array}\right)$ possibilities to choose $k$ fixed points. There remain
$N-k$ points. Let $p$ be an integer between $1$ and $N-k$. We assume
that $p$ points are directly connected to the $k$ fixed points. For
each of them, there are $k$ manners to choose one fixed point. There
are $k^p$ ways to connect directly $p$ points to $k$ fixed points.
There remains $N-k-p$ points that we must connect to the $p$ points.
There are $\sharp [\underbrace{1, \ldots, 1}_{p
~~\textrm{\scriptsize{times}}}]_{N-k}$ functions. Finally, the number of
functions with $k$ fixed points is $\displaystyle {N \choose k}
\sum\limits_{p=1}^{N-k} k^p ~\sharp [\overline{\underbrace{1,
\ldots, 1}_{p ~~\textrm{\scriptsize{times}}}}]_{N-k}$. We now prove recursively
on $N$ for every $0 \leq k \leq N$ that $\sharp
[\overline{\underbrace{1, \ldots, 1}_{k ~~\textrm{\scriptsize{times}}}}]_N =
\left(\begin{array}{c}
  N-1 \\
  N-k \\
\end{array}\right) N^{N-k}$.
We have $\sharp \overline{[1]}_1=1$. The formula is true.\\ We
suppose that $\forall k \leq N$ $\sharp [\overline{\underbrace{1,
\ldots, 1}_{k ~~\textrm{\scriptsize{times}}}}]_N = \left(\begin{array}{c}
  N-1 \\
  N-k \\
\end{array}\right) N^{N-k}$.\\
Let $X$ be a set with $N+1$ elements. We look for the functions of
$\mathcal{F}_{N+1}$ which have $k$ fixed points. Thanks to the
previous reasoning, we have
\begin{eqnarray}
\nonumber  \sharp [\overline{\underbrace{1, \ldots, 1}_{k
~~\textrm{\scriptsize{times}}}}]_{N+1} &=& \displaystyle \left(\begin{array}{c}
  N+1 \\
  k \\
\end{array}\right)
\sum\limits_{p=1}^{N+1-k} k^p ~\sharp [\overline{\underbrace{1,
\ldots, 1}_{p ~~\textrm{\scriptsize{times}}}}]_{N+1-k}.\\
 \nonumber
\sharp [\overline{\underbrace{1, \ldots, 1}_{k
~~\textrm{\scriptsize{times}}}}]_{N+1} &=& \displaystyle \left(\begin{array}{c}
  N+1 \\
  k \\
\end{array}\right)
\sum\limits_{p=1}^{N+1-k} k^p ~\sharp [\overline{\underbrace{1,
\ldots ,1}_{p ~~\textrm{\scriptsize{times}}}}]_{N-(k-1)}.
\end{eqnarray}
We use the recursion assumption.
\begin{eqnarray}
\nonumber  \sharp [\overline{\underbrace{1, \ldots, 1}_{k
~~\textrm{\scriptsize{times}}}}]_{N+1} &=& \displaystyle {N+1 \choose
k} \sum\limits_{p=1}^{N+1-k} k^p {N-k \choose p-1} (N-k+1)^{N-k+1-p}. \\
\nonumber  \sharp [\overline{\underbrace{1, \ldots, 1}_{k
~~\textrm{\scriptsize{times}}}}]_{N+1} &=& \displaystyle k {N+1 \choose
k} \sum\limits_{p=0}^{N-k}{N-k \choose p}~ k^{p} (N-k+1)^{N-k-p}. \\
\nonumber  \sharp [\overline{\underbrace{1, \ldots, 1}_{k
~~\textrm{\scriptsize{times}}}}]_{N+1}&=& \displaystyle k {N+1 \choose
k} (N+1)^{N-k} \\
\nonumber  \sharp [\overline{\underbrace{1, \ldots, 1}_{k
~~\textrm{\scriptsize{times}}}}]_{N+1} &=& \displaystyle {N \choose k-1}
(N+1)^{N-k+1} \\
\nonumber  \sharp [\overline{\underbrace{1, \ldots, 1}_{k
~~\textrm{\scriptsize{times}}}}]_{N+1} &=& \displaystyle {N \choose N-k+1}
(N+1)^{N-k+1}. \textrm{q.e.d.}
\end{eqnarray}}

\subsection{Discrete maps with $k$-cycle}

We look now for the number of functions with exactly one $k$-cycle.\newline

\stepcounter{numerothm} \textbf{{\emph{Theorem}} \thenumerothm} ~~Let $k$ be
an integer between 1 and $N$. The number of functions whose global orbit
pattern is $[k]_N$ is $\sharp [\overline{\underbrace{1, \ldots, 1}_{k
 ~~\textrm{\scriptsize{times}}}}]_{N} \times (k-1)!$.\newline

i.e. $\sharp \overline{[k]}_N = \sharp [\overline{\underbrace{1, \ldots, 1}_{%
k ~~\textrm{\scriptsize{times}}}}]_{N} \times (k-1)!$.\newline

\proof{ There are ${N \choose k}$ ways of choosing $k$ elements
among $N$. Then, there are $(k-1)!$ choices for the image of those
$k$ elements in order to constitute a $k$-cycle by $f$. We must now
count the number of ways of connecting directly or not the remaining
$N-k$ elements to the $k$-cycle. We established already this number
which is equal to $\displaystyle \sum\limits_{p=1}^{N-k} k^p ~\sharp
[\overline{\underbrace{1, \ldots, 1}_{p ~~\textrm{\scriptsize{times}}}}]_{N-k}$.
Finally, we have $\sharp \overline{[k]}_N =(k-1)! {N \choose k}
\sum\limits_{p=1}^{N-k} k^p ~\sharp [\overline{\underbrace{1,
\ldots, 1}_{p ~~\textrm{\scriptsize{times}}}}]_{N-k}$. That is, $\sharp [k]_N
=\sharp [\overline{\underbrace{1, \ldots, 1}_{k
~~\textrm{\scriptsize{times}}}}]_{N} \times (k-1)!$. q.e.d }

\subsection{Discrete maps with only two cycles}

We give the number of functions with only two cycles.\newline

\stepcounter{numerothm} \textbf{{\emph{Theorem}} \thenumerothm} ~~Let $N
\geq 2$. Let $p$ and $q$ be two non-zero integers such that $p+q \leq N$.
Then, \newline
$\sharp[\overline{p,q}]_N = \sharp [\overline{\underbrace{1, \ldots, 1}_
{p+q ~~\textrm{\scriptsize{times}}}}]_N \frac{(p+q-1)!}{q}= \frac{(N-1)!~
N^{N-(p+q)}}{(N-(p+q))!~ q} $.\newline
\newline

\proof{ We consider a function $f$ which belongs to the class
$[\overline{\underbrace{1,\ldots,1}_{p+q ~~\textrm{\scriptsize{times}}}}]_N$.
We search the number of functions constructed from $f$ whose gop is
$[p,q]_N$. From the $p$ fixed points of $f$, we construct a
$p$-cycle. Thus, there are ${{p+q-1} \choose {p-1}}$ ways to choose
$p-1$ integers among the $p+q-1$ fixed points. Counting the first
given fixed point of $f$, we have $p$ points which allow to
construct $(p-1)!$ functions with a $p$-cycle. Then there remain $q$
points which give $(q-1)!$ different functions with a $q$-cycle.
Finally, the number of functions whose gop is $[p,q]_N$ is:
${{p+q-1} \choose {p-1}} (p-1)! (q-1)!$ that is the formula
$\frac{(p+q-1)!}{q}$. }\newline

\textbf{\emph{Remark}} ~~We notice that for all $k$ non-zero integer such
that $k \leq N-1$, ~~ $\sharp \overline{[k,1]}_N = \sharp \overline{[k+1]}_N$%
.\newline

\subsection{Generalization : discrete maps with cycles of any length}

We introduce now the main theorem of the section which gives the number of
gop of discrete maps thanks to a closed formula.\newline

Given a global orbit pattern $\alpha$, the next theorem gives a formula
which gives the number of functions which belong to $\overline{\alpha}$.%
\newline

\stepcounter{numerothm} \textbf{{\emph{Theorem}} \thenumerothm} ~~Let $p
\geq 2$ be an integer. Let $[\omega_1, \ldots, \omega_p]_N$ be a gop of $%
\mathcal{G}(\mathcal{F}_N)$. Then, \newline
$\sharp \overline{[\omega_1, \ldots, \omega_p]}_N = \sharp [\overline{%
\underbrace{1, \ldots \ldots \ldots, 1}_{\omega_1 + \ldots +\omega_p ~~
\textrm{\scriptsize{times}}}}]_N \frac{(\omega_1 + \ldots +\omega_p -1)!} {\omega_p \times
(\omega_{p-1}+\omega_{p})\times \ldots \times (\omega_2 + \ldots +\omega_p)}$
.\\[.1in]
$\sharp \overline{[\omega_1, \ldots, \omega_p]}_N = \frac{(N-1)!
~N^{N-(\omega_1 + \ldots +\omega_p)}}{(N-(\omega_1 + \ldots +\omega_p))!
~\prod\limits_{k=2}^{p} (\sum\limits_{j=k}^{p} \omega_j)}$\newline
\newline

\proof{ We consider a function $f$ which belongs to $
[\overline{\underbrace{1, \ldots \ldots \ldots, 1}_{\omega_1 +
\ldots +\omega_p ~~\textrm{\scriptsize{times}}}}]_N$. We search the number of
functions constructed from $f$ whose gop is $[\omega_1, \ldots,
\omega_p]_N $. From the $\omega_1$ fixed points of $f$, we construct
a $\omega_1$-cycle. Thus, there are ${{\omega_1 + \ldots
+\omega_p-1} \choose {\omega_1-1}}$ ways to choose $\omega_1-1$
integers among the $\omega_1 + \ldots +\omega_p-1$ fixed points.
Counting the first given fixed point of $f$, we have $\omega_1$
points which allow to construct $(\omega_1-1)!$ functions with a
$\omega_1$-cycle. Then, the first fixed point of $f$ which has not
be chosen for the $\omega_1$-cycle, will belong to the
$\omega_2$-cycle. Thus, there are ${{\omega_2 + \ldots +\omega_p-1}
\choose {\omega_2-1}}$ ways to choose $\omega_2-1$ integers among
the $\omega_2 + \ldots +\omega_p-1$ fixed points. So we have
$\omega_2$ points which allow to construct $(\omega_2-1)!$ functions
with a $\omega_2$-cycle. We keep going on that way until there
remain $\omega_p$ fixed points which allow to construct
$(\omega_p-1)!$ functions with a $\omega_p$-cycle. Finally, we have
constructed : \\${{\omega_1 + \ldots +\omega_p-1} \choose
{\omega_1-1}} (\omega_1-1)! {{\omega_2 + \ldots +\omega_p-1} \choose
{\omega_2-1}} (\omega_2-1)! \times \ldots \times
{{\omega_{p-1}+\omega_p-1} \choose {\omega_{p-1}-1}}
(\omega_{p-1}-1)!(\omega_p-1)!$ functions. We simplify and obtain
the formula. }\newline

\stepcounter{numerocor} \textbf{{\emph{Corollary}} \thenumerocor} ~~Let $p$
be a non-zero integer. Let $[\omega_1, \ldots, \omega_p]_N$ be a gop of $%
\mathcal{G}(\mathcal{F}_N)$. We suppose that there exists $j$ such that $%
\omega_j \geq 2$. Let $h$ be an integer between 1 and $\omega_j-1$. Then%
\newline
$\sharp \overline{[\omega_1, \ldots, \omega_j, \ldots, \omega_p]}_N = \sharp
\overline{[\omega_1, \ldots, \omega_j-h,h, \omega_{j+1}, \ldots, \omega_p]}%
_N \times ( h + \omega_{j+1} + \ldots + \omega_p)$.\newline
\newline

\proof{ $\sharp \overline{[\omega_1, \ldots, \omega_j-h,h,
\omega_{j+1}, \ldots, \omega_p]}_N \times ( h + \omega_{j+1} +
\ldots + \omega_p) = \sharp [\overline{\underbrace{1, \ldots \ldots
\ldots, 1}_{\omega_1 + \ldots +\omega_p ~~\textrm{\scriptsize{times}}}}]_N
\\ \times \frac{(\omega_1 + \ldots +\omega_p -1)! ( h + \omega_{j+1} +
\ldots + \omega_p)}{\omega_p (\omega_{p-1}+\omega_{p}) \ldots
(\omega_{j+1}+\ldots + \omega_{p}) (h + \omega_{j+1}+\ldots +
\omega_{p}) (\omega_j + \omega_{j+1}+\ldots + \omega_{p}) \times
\ldots \times (\omega_2 +
\ldots +\omega_p)}$. \\

We simplify and we exactly obtain \\

$\sharp \overline{[\omega_1, \ldots, \omega_j-h,h, \omega_{j+1},
\ldots, \omega_p]}_N \times ( h + \omega_{j+1} + \ldots + \omega_p)
= \sharp \overline{[\omega_1,
\ldots, \omega_j, \ldots, \omega_p]}_N$. }\newline
\newline

Examples :\newline

$\sharp \overline{[2,1,3,2]}_{10} = 302,400$.\newline

$\sharp \overline{[5,2,10,8,15,2,3]}_{50} = 29, 775, 702, 147, 667, 389,
218, 762, 343, 520, 975, 006, 348,329,\newline
578, 044, 480, 000, 000, 000, 000, 000$.\newline

$\sharp \overline{[5,2,10,8,15,2,3]}_{50} \cong 2.98 \times 10^{63}$ among
the $8.88 \times 10^{84}$ functions of $\mathcal{F}_{50}$.

\section{Functions with local properties}

\label{section5}

\subsection{Functions with locally bounded range}

Since several centuries, continuity and differentiability play a
dramatic role in mathematical analysis. However these concepts are
not transposable to the functions on finite sets. Then in this
section, in order to obtain more precise results on the orbit of
dynamical systems on finite sets, we introduce subsets of
$\mathcal{F}_{N}$, whose functions have local properties such as
locally bounded range. In these subsets, the gop are found to be
fully efficient in order to describe very precisely the dynamics of
the orbits. We first consider the very simple subset
$\mathcal{L}_{1,N}$ of functions for which the difference between
$f(p)$ and $f(p+1)$ is drastically bounded. In subsection \ref{52}
we consider more sophisticated subsets.\newline

We consider the set : \newline

$\mathcal{L}_{1,N} = \{ f \in \mathcal{F}_N$ such that $\forall p, 0 \leq p
\leq N-2, |f(p)-f(p+1)| \leq 1 \}$.\newline

\subsubsection{Orbits of $\mathcal{L}_{1,N}$}

\stepcounter{numerothm} \textbf{{\emph{Theorem}} \thenumerothm} ~~If $f \in
\mathcal{L}_{1,N}$ then $f$ has only periodic orbits of order $1$ or $2$.%
\newline

\proof{We suppose that $f \in \mathcal{L}_{1,N}$ has a 3-cycle. We
denote $(a; f(a);f^2(a))$ taking $a$ the smallest value of the
3-cycle. If $a < f(a)<f^2(a)$ then there exist two non-zero integers
$e$ and $e'$ such that $f(a)=a+e$ and $f^2(a)=f(a)+e'$. Thus,
$f^2(a)-e' \leq f^3(a) \leq f^2(a)+e'$. That is $f(a)\leq a \leq
f(a)+2e'$. And finally we have the relation $a +e \leq a$ which is
impossible.\\
If $a < f^2(a)<f(a)$ then there exist two non-zero integers $e$ and
$e'$ such that $f^2(a)=a+e$ and $f(a)=f^2(a)+e'$. Thus, $f(a)-e \leq
f^3(a) \leq f(a)+e$. That is $f(a)-e\leq a \leq f(a)+e$. But
$f(a)-e=a+e'$. And finally we have the relation $a +e' \leq a$ which
is impossible.\\
We can prove in the same way that the function $f$ can't have either
3-cycle or greater order cycle than 3.}

\subsubsection{Numerical results and conjectures}

We have done numerical studies of the $\mathcal{G}(\mathcal{L}_{1,N})$ for $%
N = 1$ to $16$, using the brute force of a desktop computer (i.e. checking
every function belonging to these sets).\newline

The tables \ref{tab:g}, \ref{tab:f0}, \ref{tab:f}, \ref{tab:ff}, \ref%
{tab:fff} and \ref{tab:ffff} show the sequences for $\mathcal{L}_{1,1}$ to $%
\mathcal{L}_{1,16}$.\newline

In theses tables we display in the first column all the gop of $\mathcal{G}(%
\mathcal{L}_{1,N})$ for every value of $N$. For a given $N$, there are two
columns; the left one displays the cardinal of every existing class of gop
(- stands for non existing gop). Instead the second shows more regularity,
displaying on the row of the gop $[\underbrace{2,2,...,2}_{\mathit{k} ~
\textrm{\scriptsize{times}}}]$ the sum of the cardinals of all the classes of the gop of the form
$[\underbrace{2,2, \ldots \ldots,\underbrace{1}_{\mathit{i}^{\ th}},\ldots, 2%
}_{k+1 ~\textrm{\scriptsize{orders}}}]$ which exist.

Then we are able to formulate some statements which have not yet
been proved.

\stepcounter{numerostat} \textbf{{\emph{Statement}} \thenumerostat} ~~$%
\sharp [\overline{\underbrace{1,1, \ldots \ldots, 1}_{N-k+1 ~%
\textrm{\scriptsize{times}}}}]_{\mathcal{L}_{1,N}} =\left\{
\begin{array}{ll}
1 & {\textrm{if}}~~ k=1 \\
2 & {\textrm{if}}~~ k=2 \\
\left( \frac{4}{27}\right)(k+1) \times 3^k & {\textrm{for}}~~ 3\leq k \leq
\frac{N+1}{2}%
\end{array}
\right.$\newline
\newline

\textbf{\emph{Remark}} ~~We call $u_k= \sharp [\overline{\underbrace{1,1,
\ldots \ldots, 1}_{\mathit{N-k}+1 ~\textrm{\scriptsize{times}}}}]_{\mathcal{L}_{1,N}}
$. For $k>2$, then $u_k$ is the sequence A120926 On-line Encyclopedia of
integer Sequences : it is the number of sequences where 0 is isolated in
ternary words of length $N$ written with $\{0,1,2\}$. \newline
\newline

\stepcounter{numerostat} \textbf{{\emph{Statement}} \thenumerostat} ~~$%
\sharp [\overline{\underbrace{1,1, \ldots \ldots, 1}_{k ~
\textrm{\scriptsize{times}}}}]_{\mathcal{L}_{1,N}} =\sharp [\overline{\underbrace{1,1, \ldots \ldots, 1%
}_{k+1 ~\textrm{\scriptsize{times}}}}]_{\mathcal{L}_{1,N+1}}$ for $k \leq
\frac{N+1}{2}$.
\newline

\stepcounter{numerostat} \textbf{{\emph{Statement}} \thenumerostat} ~~$%
\sharp [\overline{\underbrace{2,2, \ldots \ldots, 2}_{k ~
\textrm{\scriptsize{times}}}}]_{\mathcal{L}_{1,N}} = \sharp [\overline{\underbrace{2,2, \ldots \ldots,
2}_{k+1 ~\textrm{\scriptsize{times}}}}]_{\mathcal{L}_{1,N+2}}$ for $k \leq
\frac{N}{2}$.
\newline

\stepcounter{numerostat} \textbf{{\emph{Statement}} \thenumerostat} ~~%
\begin{center}
$\sharp [\overline{\underbrace{2,2, \ldots \ldots, 2}_{k ~
\textrm{\scriptsize{times}}}}]_{\mathcal{L}_{1,N}} =\sharp [\overline{\underbrace{2,2,
\ldots \ldots,
2}_{k~\textrm{\scriptsize{times}}},1}]_{\mathcal{L}_{1,N+1}}$ for $2k \leq
N \leq 3k-1$.
\end{center}

\stepcounter{numerostat} \textbf{{\emph{Statement}} \thenumerostat} ~~%
\begin{center}
$\sharp [\overline{\underbrace{2,2, \ldots \ldots, 2}_{k ~
\textrm{\scriptsize{times}}}}]_{\mathcal{L}_{1,N}} = \sum\limits_{i=1}^{k+1} \sharp [\overline{%
\underbrace{2,2, \ldots \ldots,\underbrace{1}_{\mathit{i}^{\ th}},\ldots, 2}%
_{k+1 ~\textrm{\scriptsize{orders}}}}]_{\mathcal{L}_{1,N}}$ for
$2k+1 \leq N$.
\end{center}

These statements show that first the set $\mathcal{L}_{1,N}$ is an
interesting set to be considered for dynamical systems and secondly the gop
are fruitful in this study. However the set\newline

$\mathcal{L}_{2,N} = \{ f \in \mathcal{F}_N$ such that $\forall p, 0 \leq p
\leq N-2, |f(p)-f(p+1)| \leq 2 \}$\newline

\noindent is too much large to give comparable results. Then we introduce
more sophisticated sets we call sets with locally bounded range which more
or less correspond to an analogue of the discrete convolution product of the
local variation of $f$ with a compact support function $\overrightarrow{%
\alpha_t}$.

\begin{table}[tbp]
\begin{center}
\begin{tabular}{lrrrrrrrr}
\hline
g.o.p. & N=1 & N=1 & N=2 & N=2 & N=3 & N=3 & N=4 & N=4 \\
Total number &  & 1 &  & 4 &  & 17 &  & 68 \\
$[1]$ & 1 & + & 2 & + & 7 & + & 26 & + \\
$[1,1]$ & - & + & 1 & + & 4 & + & 14 & + \\
$[1,1,1]$ & - & + & - & + & 1 & + & 4 & + \\
$[1,1,1,1]$ & - & + & - & + & - & + & 1 & + \\
$[ 2]$ & - & + & 1 & 1 & 4 & 4 & 18 & 18 \\
$[2,1]$ & - & + & - & + & 1 & 1 & 3 & 4 \\
$[1,2]$ & - & + & - & + & - & + & 1 & + \\
$[2,2]$ & - & + & - & + & - & + & 1 & 1 \\ \hline
\end{tabular}%
\end{center}
\caption{Numbering functions with local properties for $f \in \mathcal{L}%
_{1,1}$, $f \in \mathcal{L}_{1,2}$, $f \in \mathcal{L}_{1,3}$, $f \in
\mathcal{L}_{1,4}$. }
\label{tab:g}
\end{table}

\begin{table}[tbp]
\begin{center}
\begin{tabular}{lrrrrrr}
\hline
g.o.p. & N=5 & N=5 & N=6 & N=6 & N=7 & N=7 \\
Total number &  & 259 &  & 950 &  & 387 \\
$[1]$ & 95 & + & 340 & + & 1,193 & + \\
$[1,1]$ & 50 & + & 174 & + & 600 & + \\
$[1,1,1]$ & 16 & + & 58 & + & 204 & + \\
$[1,1,1,1]$ & 4 & + & 16 & + & 60 & + \\
$[1,1,1,1,1]$ & 1 & + & 4 & + & 16 & + \\
$[1,1,1,1,1,1]$ & - & + & 1 & + & 4 & + \\
$[1,1,1,1,1,1,1]$ & - & + & - & + & 1 & + \\
$[ 2]$ & 70 & 70 & 264 & 264 & 952 & 952 \\
$[2,1]$ & 12 & 18 & 45 & 70 & 166 & 264 \\
$[1,2]$ & 6 & + & 25 & + & 98 & + \\
$[2,2]$ & 4 & 4 & 18 & 18 & 70 & 70 \\
$[2,2,1]$ & 1 & 1 & 4 & 4 & 17 & 18 \\
$[1,2,2]$ & - & + & - & + & 1 & + \\
$[ 2,1,2]$ & - & + & - & + & - & + \\
$[ 2,2,2]$ & - & + & 1 & 1 & 4 & 4 \\
$[2,2,2,1]$ & - & + & - & + & 1 & 1 \\ \hline
\end{tabular}%
\end{center}
\caption{Numbering functions with local properties for $f \in \mathcal{L}%
_{1,5}$, $f \in \mathcal{L}_{1,6}$, $f \in \mathcal{L}_{1,7}$. }
\label{tab:f0}
\end{table}

\begin{table}[tbp]
\small
\begin{center}
\begin{tabular}{lrrrrrr}
\hline
g.o.p. & N=8 & N=8 & N=9 & N=9 & N=10 & N=10 \\
Total number &  & 11,814 &  & 40,503 &  & 13,6946 \\
$[1]$ & 4,116 & + & 14,001 & + & 47,064 & + \\
$[ 1,1]$ & 2,038 & + & 6,852 & + & 22,806 & + \\
$[ 1,1,1]$ & 700 & + & 2,366 & + & 7,896 & + \\
$[ 1,1,1,1]$ & 214 & + & 742 & + & 2,520 & + \\
$[1,1,1,1,1]$ & 60 & + & 216 & + & 754 & + \\
$[ 1,1,1,1,1,1]$ & 16 & + & 60 & + & 216 & + \\
$[1,1,1,1,1,1,1]$ & 4 & + & 16 & + & 60 & + \\
$[1,1,1,1,1,1,1,1]$ & 1 & + & 4 & + & 16 & + \\
$[1,1,1,1,1,1,1,1,1]$ & - & + & 1 & + & 4 & + \\
$[1,1,1,1,1,1,1,1,1,1]$ & - & + & - & + & 1 & + \\
$[ 2]$ & 3,356 & 3,356 & 11,580 & 11,580 & 39,364 & 39,364 \\
$[2,1]$ & 590 & 952 & 2,062 & 3,356 & 7,072 & 11,580 \\
$[1,2]$ & 362 & + & 1,294 & + & 4,508 & + \\
$[ 2,2]$ & 264 & 264 & 952 & 952 & 3,356 & 3,356 \\
$[2,2,1]$ & 62 & 70 & 222 & 264 & 770 & 952 \\
$[1,2,2]$ & 6 & + & 28 & + & 113 & + \\
$[2,1,2]$ & 2 & + & 14 & + & 69 & + \\
$[ 2,2,2]$ & 18 & 18 & 70 & 70 & 264 & 264 \\
$[2,2,2,1]$ & 4 & 4 & 18 & 18 & 69 & 70 \\
$[1,2,2,2]$ & - & + & - & + & 1 & + \\
$[2,1,2,2]$ & - & + & - & + & - & + \\
$[2,2,1,2]$ & - & + & - & + & - & + \\
$[2,2,2,2]$ & 1 & 1 & 4 & 4 & 18 & 18 \\
$[2,2,2,2,1]$ & - & + & 1 & 1 & 4 & 4 \\
$[1,2,2,2,2]$ & - & + & - & + & - & + \\
$[2,1,2,2,2]$ & - & + & - & + & - & + \\
$[2,2,1,2,2]$ & - & + & - & + & - & + \\
$[2,2,2,1,2]$ & - & + & - & + & - & + \\
$[2,2,2,2,2]$ & - & + & - & + & 1 & 1 \\ \hline
\end{tabular}%
\end{center}
\caption{Numbering functions with local properties for $f \in \mathcal{L}%
_{1,8}$, $f \in \mathcal{L}_{1,9}$, $f \in \mathcal{L}_{1,10}$. }
\label{tab:f}
\end{table}

\begin{table}[tbp]
\small 
\begin{center}
\begin{tabular}{lrrrrrr}
\hline
g.o.p. & N=11 & N=11 & N=12 & N=12 & N=13 & N=13 \\
Total number &  & 457,795 &  & 1,515,926 &  & 4,979,777 \\
$[1]$ & 156,629 & + & 516,844 & + & 1,693,073 & + \\
$[1,1]$ & 75,292 & + & 246,762 & + & 803,706 & + \\
$[1,1,1]$ & 26,098 & + & 85,556 & + & 278,580 & + \\
$[1,1,1,1]$ & 8,434 & + & 27,904 & + & 91,488 & + \\
$[1,1,1,1,1]$ & 2,756 & + & 8,658 & + & 28,738 & + \\
$[1,1,1,1,1,1]$ & 756 & + & 2,590 & + & 8,730 & + \\
$[1,1,1,1,1,1,1]$ & 216 & + & 756 & + & 2,592 & + \\
$[1,1,1,1,1,1,1,1]$ & 60 & + & 216 & + & 756 & + \\
$[1,1,1,1,1,1,1,1,1]$ & 16 & + & 60 & + & 216 & + \\
$[1,1,1,1,1,1,1,1,1,1]$ & 4 & + & 16 & + & 60 & + \\
$[1,1,1,1,1,1,1,1,1,1,1]$ & 1 & + & 4 & + & 16 & + \\
$[1,1,1,1,1,1,1,1,1,1,1,1]$ & - & + & 1 & + & 4 & + \\
$[1,1,1,1,1,1,1,1,1,1,1,1,1]$ & - & + & - & + & 1 & + \\
$[2]$ & 132,104 & 132,104 & 438,846 & 438,846 & 1,445,258 & 1,445,258 \\
$[2,1]$ & 23,941 & 39,364 & 80,108 & 132,104 & 265,548 & 438,846 \\
$[1,2]$ & 15,423 & + & 51,996 & + & 173,298 & + \\
$[2,2]$ & 11,580 & 11,580 & 39,364 & 39,364 & 132,104 & 132,104 \\
$[2,2,1]$ & 2,634 & 3,356 & 8,883 & 11,580 & 29,659 & 39,364 \\
$[1,2,2]$ & 429 & + & 1,555 & + & 5,478 & + \\
$[2,1,2]$ & 293 & + & 1,142 & + & 4,227 & + \\
$[2,2,2]$ & 952 & 952 & 3,356 & 3,356 & 11,580 & 11,580 \\
$[2,2,2,1]$ & 255 & 264 & 899 & 952 & 3,098 & 3,356 \\
$[ 1,2,2,2]$ & 7 & + & 35 & + & 152 & + \\
$[ 2,1,2,2]$ & 2 & + & 16 & + & 86 & + \\
$[2,2,1,2]$ & - & + & 2 & + & 20 & + \\
$[2,2,2,2]$ & 70 & 70 & 264 & 264 & 952 & 952 \\
$[ 2,2,2,2,1]$ & 18 & 18 & 70 & 70 & 263 & 264 \\
$[ 1,2,2,2,2]$ & - & + & - & + & 1 & + \\
$[2,1,2,2,2]$ & - & + & - & + & - & + \\
$[2,2,1,2,2]$ & - & + & - & + & - & + \\
$[2,2,2,1,2]$ & - & + & - & + & - & + \\
$[2,2,2,2,2]$ & 4 & 4 & 18 & 18 & 70 & 70 \\
$[2,2,2,2,2,1]$ & 1 & 1 & 4 & 4 & 18 & 18 \\
$[1,2,2,2,2,2]$ & - & + & - & + & - & + \\
$[2,2,2,2,2,2]$ & - & + & 1 & 1 & 4 & 4 \\
$[2,2,2,2,2,2,1]$ & - & + & - & + & 1 & 1 \\ \hline
\end{tabular}%
\end{center}
\caption{Numbering functions with local properties for $f \in \mathcal{L}%
_{1,11}$, $f \in \mathcal{L}_{1,12}$, $f \in \mathcal{L}_{1,13}$.}
\label{tab:ff}
\end{table}

\begin{table}[tbp]
\small 
\begin{center}
\begin{tabular}{lrrrr}
\hline
g.o.p. & N=14 & N=14 & N=15 & N=15 \\
Total number &  & 16,246,924 &  & 52,694,573 \\
$[1]$ & 5,511,218 & + & 17,841,247 & + \\
$[1,1]$ & 2,603,258 & + & 8,391,360 & + \\
$[ 1,1,1]$ & 901,802 & + & 2,904,592 & + \\
$[1,1,1,1]$ & 297,728 & + & 962,888 & + \\
$[1,1,1,1,1]$ & 94,440 & + & 307,848 & + \\
$[ 1,1,1,1,1,1]$ & 29,050 & + & 95,676 & + \\
$[1,1,1,1,1,1,1]$ & 8,746 & + & 29,140 & + \\
$[1,1,1,1,1,1,1,1]$ & 2,592 & + & 8,748 & + \\
$[ 1,1,1,1,1,1,1,1,1]$ & 756 & + & 2,592 & + \\
$[1,1,1,1,1,1,1,1,1,1]$ & 216 & + & 756 & + \\
$[1,1,1,1,1,1,1,1,1,1,1]$ & 60 & + & 216 & + \\
$[ 1,1,1,1,1,1,1,1,1,1,1,1]$ & 16 & + & 60 & + \\
$[ 1,1,1,1,1,1,1,1,1,1,1,1,1]$ & 4 & + & 16 & + \\
$[ 1,1,1,1,1,1,1,1,1,1,1,1,1,1]$ & 1 & + & 4 & + \\
$[1,1,1,1,1,1,1,1,1,1,1,1,1,1,1]$ & - & + & 1 & + \\
$[2]$ & 4,725,220 & 4,725,220 & 15,352,392 & 15,352,392 \\
$[2,1]$ & 873,149 & 1,445,258 & 2,851,350 & + \\
$[1,2]$ & 572,109 & + & 1,873,870 & + \\
$[2,2]$ & 438,846 & 438,846 & 1,445,258 & 1,445,258 \\
$[2,2,1]$ & 98,135 & 132,104 & 322,310 & 438,846 \\
$[1,2,2]$ & 18,873 & + & 63,967 & + \\
$[2,1,2]$ & 15,096 & + & 52,569 & + \\
$[2,2,2]$ & 39,364 & 39,364 & 132,104 & 132,104 \\
$[2,2,2,1]$ & 10,460 & 11,580 & 34,845 & 39,364 \\
$[1,2,2,2]$ & 605 & + & 2,282 & + \\
$[2,1,2,2]$ & 389 & + & 1,596 & + \\
$[ 2,2,1,2]$ & 126 & + & 641 & + \\
$[2,2,2,2]$ & 3,356 & 3,356 & 11,580 & 11,580 \\
$[2,2,2,2,1]$ & 942 & 952 & 3,292 & 3,356 \\
$[1,2,2,2,2]$ & 8 & + & 44 & + \\
$[2,1,2,2,2]$ & 2 & + & 18 & + \\
$[2,2,1,2,2]$ & - & + & 2 & + \\
$[2,2,2,1,2]$ & - & + & - & + \\
$[2,2,2,2,2]$ & 264 & 264 & 952 & 952 \\
$[2,2,2,2,2,1]$ & 70 & 70 & 264 & 264 \\
$[1,2,2,2,2,2]$ & - & + & - & + \\
$[2,2,2,2,2,2]$ & 18 & 18 & 70 & 70 \\
$[2,2,2,2,2,2,1]$ & 4 & 4 & 18 & 18 \\
$[2,2,2,2,2,2,2]$ & 1 & 1 & 4 & 4 \\
$[2,2,2,2,2,2,2,1]$ & - & + & 1 & 1 \\ \hline
\end{tabular}%
\end{center}
\caption{Numbering functions with local properties for $f \in \mathcal{L}%
_{1,14}$, $f \in \mathcal{L}_{1,15}$. }
\label{tab:fff}
\end{table}

\begin{table}[tbp]
\small
\begin{center}
\begin{tabular}{lrr}
\hline
g.o.p. & N=16 & N=16 \\
Total number &  & 170,028,792 \\
$[1]$ & 57,477,542 & + \\
$[1,1]$ & 26,932,398 & + \\
$[ 1,1,1]$ & 9,314,088 & + \\
$[1,1,1,1]$ & 3,097,650 & + \\
$[1,1,1,1,1]$ & 996,764 & + \\
$[ 1,1,1,1,1,1]$ & 312,456 & + \\
$[1,1,1,1,1,1,1]$ & 96,096 & + \\
$[1,1,1,1,1,1,1,1]$ & 29,158 & + \\
$[ 1,1,1,1,1,1,1,1,1]$ & 8,748 & + \\
$[1,1,1,1,1,1,1,1,1,1]$ & 2,592 & + \\
$[1,1,1,1,1,1,1,1,1,1,1]$ & 756 & + \\
$[ 1,1,1,1,1,1,1,1,1,1,1,1]$ & 216 & + \\
$[ 1,1,1,1,1,1,1,1,1,1,1,1,1]$ & 60 & + \\
$[ 1,1,1,1,1,1,1,1,1,1,1,1,1,1]$ & 16 & + \\
$[1,1,1,1,1,1,1,1,1,1,1,1,1,1,1]$ & 4 & + \\
$[1,1,1,1,1,1,1,1,1,1,1,1,1,1,1,1]$ & 1 & + \\
$[2]$ & 49,610,818 & 49,610,818 \\
$[2,1]$ & 9,255,822 & 15,352,392 \\
$[1,2]$ & 6,096,570 & + \\
$[2,2]$ & 4,725,220 & 4,725,220 \\
$[2,2,1]$ & 1,051,686 & 1,445,258 \\
$[1,2,2]$ & 213,975 & + \\
$[2,1,2]$ & 179,597 & + \\
$[2,2,2]$ & 438,846 & 438,846 \\
$[2,2,2,1]$ & 114,798 & 132,104 \\
$[1,2,2,2]$ & 8,284 & + \\
$[2,1,2,2]$ & 6,146 & + \\
$[ 2,2,1,2]$ & 2,876 & + \\
$[2,2,2,2]$ & 39,364 & 39,364 \\
$[2,2,2,2,1]$ & 11,246 & 11,580 \\
$[1,2,2,2,2]$ & 204 & + \\
$[2,1,2,2,2]$ & 106 & + \\
$[2,2,1,2,2]$ & 22 & + \\
$[2,2,2,1,2]$ & 2 & + \\
$[2,2,2,2,2]$ & 3,356 & 3,356 \\
$[2,2,2,2,2,1]$ & 951 & 952 \\
$[1,2,2,2,2,2]$ & 1 & + \\
$[2,2,2,2,2,2]$ & 264 & 264 \\
$[2,2,2,2,2,2,1]$ & 70 & 70 \\
$[2,2,2,2,2,2,2]$ & 18 & 18 \\
$[2,2,2,2,2,2,2,1]$ & 4 & 4 \\
$[2,2,2,2,2,2,2,2]$ & 1 & 1 \\ \hline
\end{tabular}%
\end{center}
\caption{Numbering functions with local properties for $f \in \mathcal{L}%
_{1,16}$. }
\label{tab:ffff}
\end{table}

\begin{table}[tbp]
\small
\begin{center}
\begin{tabular}{rrrrr}
\hline
q & maximal period & modulus & gop number & functions number \\
35 & 2 & 2 & 3 & 9,992 \\
41 & 2 & 3 & 6 & 21,764 \\
48 & 3 & 3 & 7 & 63,408 \\
51 & 3 & 4 & 9 & 122,316 \\
54 & 4 & 4 & 15 & 258,910 \\
58 & 4 & 5 & 19 & 497,106 \\
60 & 4 & 5 & 25 & 696,586 \\
61 & 4 & 10 & 37 & 818,000 \\
62 & 4 & 10 & 44 & 921,698 \\
63 & 4 & 10 & 46 & 1,022,184 \\
68 & 5 & 10 & 50 & 1,604,518 \\
70 & 6 & 10 & 60 & 1,837,088 \\
72 & 6 & 10 & 61 & 2,124,974 \\
73 & 6 & 10 & 88 & 2,352,560 \\
76 & 6 & 10 & 98 & 3,514,608 \\
77 & 6 & 10 & 99 & 4,001,306 \\
81 & 6 & 10 & 100 & 6,499,244 \\
82 & 6 & 10 & 104 & 7,230,576 \\
83 & 7 & 10 & 109 & 8,113,212 \\
84 & 8 & 10 & 117 & 9,126,054 \\
85 & 8 & 10 & 130 & 10,184,542 \\
86 & 8 & 10 & 131 & 11,244,702 \\
87 & 8 & 10 & 145 & 12,311,866 \\
88 & 8 & 10 & 161 & 13,485,506 \\
89 & 8 & 10 & 175 & 14,692,658 \\
90 & 8 & 10 & 176 & 15,984,782 \\
92 & 8 & 10 & 182 & 18,775,284 \\
93 & 8 & 10 & 188 & 20,252,084 \\
94 & 8 & 10 & 193 & 21,640,666 \\
95 & 8 & 10 & 195 & 23,021,112 \\
96 & 8 & 10 & 242 & 24,479,312 \\
97 & 8 & 10 & 298 & 26,163,582 \\
98 & 8 & 10 & 335 & 28,285,274 \\
99 & 9 & 10 & 377 & 30,861,396 \\
100 & 10 & 10 & 463 & 34,086,310 \\
101 & 10 & 10 & 484 & 37,553,504 \\ \hline
\end{tabular}%
\end{center}
\caption{Numerical study of the set $\mathcal{L}_{\vec{\protect\alpha_t},q,N}
$ for $N=10$, $t=5$, $\protect\alpha_1=20$, $\protect\alpha_2=10$, $\protect%
\alpha_3=5$, $\protect\alpha_4=3$ and $\protect\alpha_5=1$, for $q=35, ...,
101$}
\label{tab:ggg}
\end{table}

\normalsize 

\subsection{Orbits and gop of locally bounded range function sets}

\label{52}

Consider now the set :

$\mathcal{L}_{\overrightarrow{\alpha_t},q,N} = \{ f \in
\mathcal{F}_N$ such that $\forall p, 0 \leq p \leq N-r-1,
\sum\limits_{r=1}^{r=t} \alpha_r |f(p)-f(p+r)| \leq q \} \bigcap \{
f \in \mathcal{F}_N$ such that $\forall p, t \leq p \leq N-1,
\sum\limits_{r=1}^{r=t} \alpha_r |f(p)-f(p-r)| \leq q \}$ for the
vector $\overrightarrow{\alpha_t}=(\alpha_1, \alpha_2, \ldots,
\alpha_t) \in \mathbb{N}^t$, for $q \in \mathbb{N}$.\\

The functions belonging to these sets show a kind of "rigidity": the less is
$q$, the more "rigid" is the function. Furthermore, the maximal length of a
periodic orbit increases with $q$, and so the number of gop $\sharp \mathcal{%
G}(\mathcal{L}_{\overrightarrow{\alpha_t},q,N})$ and the maximal modulus of
the gop. As an example, we explore numerically the case: $N=10$, $t=5$, $%
\alpha_1=20$, $\alpha_2=10$, $\alpha_3=5$, $\alpha_4=3$ and
$\alpha_5=1$, for $q=35, \ldots, 101$. The results are displayed in
table \ref{tab:ggg}. In this table "modulus" means the maximal
modulus of the gop belonging to this set for the corresponding value
of $q$ in the row, "gop number" stands for $\sharp
\mathcal{G}(\mathcal{L}_{\overrightarrow{\alpha_t},q,N})$ and
"functions number" for $\sharp
\mathcal{L}_{\overrightarrow{\alpha_t},q,N}$. One can point out that
for the particular function $\overrightarrow{\alpha_t} $ of the
example; it is possible to find $10$ intervals $I_1, I_2, \ldots,
I_{10} \subset \mathbb{N}$ such that if $q \in I_r$ then there is no
periodic orbit whose period is strictly greater than $r$, (e.g.,
$I_6= [\![ 70,82 ]\!]$). Furthermore it is possible to split these
intervals into
subintervals $I_{r,s}$ in which $\sharp \mathcal{G}(\mathcal{L}_{%
\overrightarrow{\alpha_t},q,N})$ is constant when $q$ thumbs $I_{r,s}$. This
is not the case for $\sharp \mathcal{L}_{\overrightarrow{\alpha_t},q,N}$.

\section{Conclusion}

\label{section6} The behaviour of a discrete dynamical system associated to
a function on finite ordered set $X$ is not easily predictable. Such a
system can only exhibit periodic orbits. All the orbits can be explicitly
computed; all together they form a global orbit pattern. We formalise such a
gop as the ordered set of periods when the initial value thumbs $X$ in the
increasing order. We are able to predict, using closed formulas, the number
of gop for the set $\mathcal{F}_N$ of all the functions on $X$. We introduce
special subsets of $\mathcal{F}_N$ in order to understand more precisely the
behaviour of the dynamical system. We explore by computer experiments theses
sets, they show interesting patterns of gop. Further study is needed to
understand the behaviour of dynamical systems associated to functions
belonging to these sets.

\end{document}